\documentclass[12pt]{amsart}

\newsymbol\pp 1275
\newsymbol\twoheadrightarrow 1310

\newtheorem{theorem}{Theorem}[section]
\newtheorem{proposition}[theorem]{Proposition}

\newtheorem{lemma}[theorem]{Lemma}
\newtheorem{corollary}[theorem]{Corollary}

\theoremstyle{definition}

\newcommand{\Q}{{\mathbb Q}}

\newcommand{\Z}{{\mathbb Z}}

\newcommand{\mcL}{ {\mathcal L}}
\newcommand{\mcI}{ {\mathcal I}}
\newcommand{\mcA}{ {\mathcal A}}

\DeclareMathOperator{\Mod}{Mod} \DeclareMathOperator{\Aut}{Aut} \DeclareMathOperator{\Hom}{Hom}
\DeclareMathOperator{\End}{End} \DeclareMathOperator{\Sp}{Sp} \DeclareMathOperator{\krn}{ker}

\usepackage{color}      
\usepackage{epsfig}
\usepackage{amsfonts}
\usepackage{psfrag}

\title{Pseudo-Anosov homeomorphisms and the lower central series of a surface group}
\author{Justin Malestein}

\begin{document}
 \maketitle

Version: February 8th, 2007

\section{Introduction}
\noindent Denote by $\Mod(S)$ the mapping class group of a compact, oriented surface $S = S_{g,1}$ of genus $g \geq 2$ with one
boundary component; i.e. $\Mod(S)$
 is the group of homeomorphisms of $S$ fixing $\partial S$ pointwise up to isotopies fixing $\partial S$ pointwise.
  A basic question to
contemplate is: what topological or dynamical data of a mapping class can be extracted from various kinds of
algebraic data?  Since pseudo-Anosovs are the more complex mapping classes topologically and dynamically, we
would like to know if a given mapping class is pseudo-Anosov; i.e. it has a representative homeomorphism
which leaves invariant a pair of transverse measured foliations.

One kind of algebraic data is the action of a mapping class on $\Gamma := \pi_1(S, *)$ and its
various quotients. Specifically, consider the sequence of $k$-step nilpotent quotients
$N_k := \Gamma
/ \Gamma_{k+1}$ where $\{\Gamma_k\}$ is the lower central series of $\Gamma$ defined inductively by
$$\Gamma_1 = \Gamma  \;\;\;\;\;\;\;\;  \Gamma_k = [\Gamma, \Gamma_{k-1}] \;\;  \textrm{for} \;\; k > 1 $$
Since elements of $\Mod(S)$ fix $\partial S$ pointwise and we choose the basepoint $* \in \partial S$, we
obtain a representation $\Mod(S) \to \Aut(\Gamma)$, and
furthermore since each $\Gamma_k$ is characteristic, we obtain a representation for each $k$:
$$ \rho_k: \Mod(S) \to \Aut(\Gamma / \Gamma_{k+1})$$
One natural question to ask is: given {\it only} the datum of $\rho_k(f)$ for $f \in \Mod(S)$, can we determine
if the mapping class is pseudo-Anosov or not?  If the mapping class is determined to be pseudo-Anosov, can we
detect the dilatation?  This paper is one step in answering the first question.

For $k \geq 1$, we define the $k$th Torelli group to be $\mcI_k(S) := \krn(\rho_k)$ (and so with our indexing, which is
different from some other authors, the
classical Torelli group is $\mcI_1(S)$).  To each $f \in \mcI_k = \mcI_k(S)$, we associate an invariant $\Psi_k(f) \in
\End(H_1(S, \Z))$
which is constructed from $\rho_{k+1}(f)$.  We will prove the following:

\begin{theorem}[Criterion for pseudo-Anosovs]
 Let $f \in \mcI_k$.  If the characteristic polynomial of $\Psi_k(f)$ is
irreducible in $\Z[x]$, then $f$ is pseudo-Anosov.
\end{theorem}

\noindent This follows immediately from the following theorem which we prove in Section
\ref{section:proofofcriterion}.  For the remainder of this paper, we let $H := H_1(S, \Z)$.
\begin{theorem}
\label{theorem:maintheorem} Let $f \in \mcI_k$.  If the characteristic polynomial $\chi(\Psi_k(f))$ of $\Psi_k(f) \in
\End(H)$ has no (nontrivial) even degree or degree 1 factors over $\Z$, then $f$ is pseudo-Anosov.
\end{theorem}

\noindent Since $\Psi_k$ uses only the data of $\rho_{k+1}(f)$ and $\krn(\rho_{k+1}) = \mcI_{k+1}$, we obtain the
following corollary:
\begin{corollary}
 If $f \in \mcI_k$ satisfies the hypothesis of Theorem \ref{theorem:maintheorem}, then the whole coset $f \mcI_{k+1}$
is pseudo-Anosov.
\end{corollary}

Note that the data of $\rho_1$ is not used in Theorem \ref{theorem:maintheorem}.  Since $\Gamma / \Gamma_2 =
H $, the homomorphism $\rho_1$ is the standard representation into $\Aut(H)$ with image isomorphic to
the integral symplectic group $\Sp(2g, \Z)$.  It is not too difficult to find a criterion on $\rho_1(f)$ for $f$ to
be pseudo-Anosov, and in fact, Casson and Bleiler give such a criterion in Lemma 5.1 of \cite{CB}.
  Casson--Bleiler show that if the characteristic polynomial, $\chi(\rho_1(f))$, is irreducible over $\Z$, has no
roots of unity as eigenvalues, and is not $g(t^n)$ for any $n > 1$ and $g \in \Z[x]$, then $f$ is pseudo-Anosov.

The Casson--Bleiler criterion is well-known and has been around for many years.  It is unfortunately unable to detect
pseudo-Anosovs in any of the $\mcI_k$ simply because $\mcI_k \subseteq \krn(\rho_1)$.  (This is not to imply that the
Casson-Bleiler criterion can detect all pseudo-Anosovs which act non-trivially on $H$.)  In this sense, this paper
is filling in the gap that Casson-Bleiler left (for the case of a surface with one boundary component). \\

\noindent \textbf{Remark:}  It is well-known that $\mcI_1$ has pseudo-Anosov elements thanks to criteria
of Thurston, Penner, and others \cite{T} \cite{P} \cite{BH}.  However, their methods of finding pseudo-Anosovs are all
topological as opposed to algebraic in nature.  Furthermore, their criteria require the specification of a particular
mapping class and thus are not well-suited to dealing with the information of $\rho_k(f)
\in \Aut(\Gamma_1 / \Gamma_{k+1})$ which only specifies a coset of $\mcI_k$.   Both Thurston's criterion and Penner's
criterion require that a
mapping class be described in terms of twists about two multi-curves.  In \cite{BH}, Bestvina--Handel describe an
algorithm using train tracks that can determine whether any single mapping class is pseudo-Anosov or not.  In fact,
this algorithm has been implemented in a computer program by Peter Brinkmann (\cite{Br}).\\

Let us now outline the contents of the paper.  In Section  \ref{section:basicfacts}, we recall some basic properties
of the series $\{\Gamma_k\}$. We then define for $f \in \mcI_k$ the invariant $\Psi_k(f) \in \End(H)$. ($\Psi_k(f)$
is in general non-trivial which might be rather suprising given that $\rho_1(f) \in \Aut(H)$ is necessarily trivial.)
To define $\Psi_k$, we need two ingredients: the Johnson homomorphism $\tau$ and contractions
$$\Phi_{2k} : \Gamma_{2k+1} / \Gamma_{2k+2} \to H$$
Defining $\Phi_{2k}$ requires a bit of work and is described in Section \ref{section:contractionsection}.  In Section
\ref{section:definingtau}, we will recall the definition of the Johnson homomorphism $\tau$ which we describe here as
follows:
$$\tau: \mcI_k(S) \to \Hom\left(\bigoplus_{m=1}^\infty \Gamma_{m} /
\Gamma_{m+1}, \bigoplus_{m=k+1}^\infty \Gamma_{m} / \Gamma_{m+1}\right)$$  We will denote the image of $f$ under $\tau$
as $\tau_f$.  We define $\Psi_k$ as follows:
$$\Psi_k(f) := \left\{ \begin{array}{ll} \Phi_k \circ (\tau_f|_H)  & \textrm{$k$ even} \\
                                         \Phi_{2k} \circ (\tau_f^2|_H)  & \textrm{$k$ odd}  \end{array} \right. \in
                                     \End(H)$$
Note that the map $\Psi_k$ is a homomorphism for $k$ even but not necessarily for $k$ odd.

In Section \ref{section:proofofcriterion}, we prove Theorem \ref{theorem:maintheorem}.  The general idea of the proof
of Theorem \ref{theorem:maintheorem} is to use the Nielsen-Thurston classification which states that a mapping class
is pseudo-Anosov if and only if it is neither reducible nor of finite order.  Recall that $f$ is {\em reducible} if
$f$ fixes the isotopy class of an essential 1-dimensional submanifold where {\em essential} means that each component
is neither null-homotopic nor homotopic to a boundary component.  Since $\mcI_1$ is torsion-free, the classification
 reduces to: $f$
is pseudo-Anosov if and only if it is irreducible. We then show that reducibility of $f$ implies that
$\chi(\Psi_k(f))$ has a linear or even degree factor by using the fact that a certain subgroup of $\pi_1(S)$ is
invariant under $f_* \in \Aut(\pi_1(S))$.

For any particular $f \in \mcI_k$, the invariant $\Psi_k(f)$ is explicitly computable, provided one can compute $\tau_f$.
  In Section \ref{section:example}, we show some mapping classes satisfy the hypothesis
of Theorem \ref{theorem:maintheorem} by computing $\Psi_k(f)$ directly.  Nevertheless, at present the author has not
found whole families of pseudo-Anosovs ranging over either $g$ or $k$ which satisfy the hypothesis of Theorem
\ref{theorem:maintheorem}.  Additionally, in section \ref{section:example} we compare Theorem
\ref{theorem:maintheorem} to the Thurston/Penner criteria. \\

\noindent \textbf{Remark:}  We choose to work with a surface with a boundary component as opposed to a closed surface
to simplify things technically.  The fundamental group of a surface with boundary is a free group.  As we shall see
in Section \ref{section:basicfacts}, this will further imply that the Lie algebra associated to the $\{\Gamma_k\}$
 is a free Lie algebra.  While the author suspects that one may obtain a criterion for closed surfaces
from this criterion, he has not done so at present. \\

\noindent \textbf{Acknowledgements.}  The author would like to thank Dan Margalit, Nathan Broaddus, Ian Biringer,
Juan Souto, Matthew Day and Asaf Hadari for their helpful comments.  He would also like to thank Andy Putman for help
during the research stage.  He would like to especially thank Benson Farb for extensive comments, posing the
question, continuous help and inspiration.

\section{Basic facts about the lower central series} \label{section:basicfacts} \noindent  For the reader's convenience,
we recall basic facts about central filtrations of a group.  Suppose $$ G = G_1 \supset G_2 \supset G_3 \dots$$
is a filtration of $G$ by normal subgroups.  We call $G$ a {\em central filtration} if $[G_k, G_l] \subseteq
G_{k+l}$.  We recount the following folklore result.
\begin{theorem} \label{theorem:liealg}
Let $\{G_i\}$ be a central filtration of $G$ by normal subgroups.  Then, the following hold:
\begin{enumerate}
\item The function $G_k \times G_l \to G_{k+l}$ given by $(x, y) \mapsto xyx^{-1}y^{-1}$ induces a
well-defined map $$G_k / G_{k+1} \times G_l / G_{l+1} \to G_{k+l} / G_{k+l+1}$$
\item Using the pairing from (1) as a bracket which we denote by $[\,\,,\, ]$, we obtain a graded $\Z$-Lie algebra:
$$L := \bigoplus_k G_k / G_{k+1}$$
\end{enumerate}
\end{theorem}
For an explanation and proof see Sections 3.1 and 4.5 of \cite{BL}.  Also, we recall for the reader that the lower
central series is a central filtration (see 4.4 of \cite{BL}).

 Recall that the fundamental group of a surface with
boundary is a free group.  The Lie algebra associated to a free group's lower central series is special as described
in the following theorem which is a rephrasing of Theorem 5.12 of \cite{MKS}.
\begin{theorem} \label{theorem:freeliealg}
Let $G$ be a free group with generators $a_1, \dots, a_n$ and lower central series $G_1 \supset G_2 \supset \dots$.
Then the (graded) $\Z$-Lie algebra $$L: = (\bigoplus_k G_k / G_{k+1},\,\,\, [\,\,,\,]) $$ is a free $\Z$-Lie
algebra. $L$ has as its generating set $\{a_1, \dots, a_n\}$ viewed as a subset of $G_1 / G_2$.
\end{theorem}
The definition of {\em free Lie algebra} is exactly what one expects: given a $\Z$-Lie algebra $L'$ and elements $x_1,
\dots, x_n \in L'$, there exists a unique Lie algebra homomorphism $h: L \to L'$ such that $h(a_i) = x_i$. The free
Lie algebra in general is fairly complicated.  Even computing the rank of $G_k / G_{k+1}$ for arbitrary $k$ is
nontrivial.  Thankfully, free Lie algebras embed in simpler Lie algebras.

A {\em free associative $\Z$-algebra} $A$ with generators $b_1, \dots, b_n$ is a noncommutative ring with the
universal property that given a $\Z$-algebra $A'$ and elements $x_1, \dots, x_n \in A'$ there is a unique
homomorphism $h: A \to A'$ such that $h(b_i) = x_i$.  More concretely, $A$ is (canonically isomorphic to)
the noncommutative polynomial ring in $n$ variables over $\Z$.  However, viewing $A$ as a polynomial ring is not
particularly convenient for the purposes of this paper.  If we let $M := \Z^n$, then $A$ is isomorphic to the
tensor algebra $ \bigoplus_{k=0}^{\infty} M^{\otimes k}$ where $M^{\otimes 0} := \Z$.  The algebra $A$ has a
canonical Lie bracket: $[x, y] := x \otimes y - y \otimes x$.  Thus, we
have a canonical Lie homomorphism $ \mcL  \to   A $ defined by $a_i \mapsto b_i$.  From Corollary 0.3 and Theorem
0.5 of \cite{R}, we obtain the following.
\begin{theorem} \label{theorem:freeassocalg}
If $L$ is a free $\Z$-Lie algebra with generators $a_1, \dots, a_n$ and $A$ is a free associative algebra over $\Z$
with generators $b_1, \dots, b_n$, then the canonical Lie homomorphism induced by $a_i \mapsto b_i$ is injective.
\end{theorem}
\noindent Moreover, it is not hard to check that the map $L \to A$ respects the grading.

Now, let us apply Theorems \ref{theorem:freeliealg} and \ref{theorem:freeassocalg} to the group $\Gamma := \pi_1(S)$
with (free) generators $a_1, \dots, a_{2g}$.  Let $\mcL$
 be the graded Lie algebra associated to $\{\Gamma_k\}$.  Let $\mcA$ be the tensor algebra $ \bigoplus_{k=0}^\infty
H^{\otimes k}$ where $H^{\otimes 0}:= \Z$.  Since $H \cong \Z^{2g}$, the algebra $\mcA$ is a free associative algebra
To simplify notation, let us define $\mcL_k := \Gamma_k/\Gamma_{k+1}$.  Recall that $\mcA \cong
\bigoplus_{k=0}^{\infty} M^{\otimes k}$ where $M = \Z^{2g}$.  We have defined the $a_i$ as elements
of $\pi_1(S)$, but we can also consider the equivalence class of $a_i$ in $\Gamma_1 / \Gamma_2 \subset \mcL$ or in $H
= H^{\otimes 1} \subset \mcA$.  Thus, we obtain a natural, injective map $\mcL \to \mcA$ defined by sending ``$a_i$'' to
``$a_i$''.

The mapping class group has a natural action on $\mcL$ by considering $$\mcL = \bigoplus_{k=1}^\infty \Gamma_k /
\Gamma_{k+1}$$ as a direct sum of representations $\Mod(S) \to \Aut(\Gamma_k / \Gamma_{k+1})$.  We obtain an action
on $$\mcA = \bigoplus_{k=0}^\infty  H^{\otimes k}$$ from the action on $H$.  It is not hard to check that the map
$\mcL \to \mcA$ respects this action.  Since the $\Mod(S)$-action on $\mcA$ is induced by the action on $H$, it
factors through to an $\Sp(2g, \Z)$-action and so the $\Mod(S)$-action on $\mcL$ factors through $\Sp(2g, \Z)$ also
(This can also be proven directly.).

\section{The Johnson Homomorphisms} \label{section:definingtau}
All of the results in this section are the work of Johnson, Morita, Hain and others.  Recall that
$$\mcI_k := \krn(\Mod(S) \to \Aut(\Gamma_1 / \Gamma_{k+1}))$$ and $H = H_1(S)$.  A preliminary version of the Johnson
homomorphism is a map: $$\tau: \mcI_k \to \Hom(H, \Gamma_{k+1} / \Gamma_{k+2})$$ for each $k$.  Note that the image
of $f$ under $\tau$ will be denoted $\tau_f$ as is standard.  We define the preliminary version as follows.  Let $f \in
\mcI_k$. Since $f_*$ acts trivially on $\Gamma_1 / \Gamma_{k+1}$, we obtain a well-defined map of sets
$$\begin{array}{rcl} t_f: \Gamma_1 / \Gamma_{k+2} &  \to & \Gamma_{k+1} / \Gamma_{k+2}\\  x & \mapsto & f_*(x) x^{-1}
\end{array}$$  The following result is one part of Proposition 2.3 in \cite{M3}.

\begin{proposition}[Johnson, Morita] \label{proposition:basictau}
The set map $t_f: \Gamma_1 / \Gamma_{k+2} \to \Gamma_{k+1} / \Gamma_{k+2}$ induces a well-defined homomorphism $H \to
\Gamma_{k+1} / \Gamma_{k+2}$ which is $\tau_f$.  Moreover, $\tau$ is a homomorphism.
\end{proposition}

\begin{proof}
By the very definition of the lower central series, $\Gamma_{k+1} / \Gamma_{k+2}$ is in the center of $\Gamma_1 /
\Gamma_{k+2}$.
  Thus, $$f_*(xy)(xy)^{-1} = f_*(xy)y^{-1}x^{-1} = f_*(x)(f_*(y)y^{-1})x^{-1} = f_*(x)x^{-1}(f_*(y)y^{-1})$$ and so
$t_f$ is in fact a homomorphism. As $\Gamma_{k+1} / \Gamma_{k+2}$ is abelian, this homomorphism factors through the
abelianization of $\Gamma_1 / \Gamma_{k+1}$ which is $\Gamma_1 / [\Gamma_1, \Gamma_1] = \Gamma_1 / \Gamma_2 = H$.
Hence, we obtain a homomorphism $H \to \Gamma_{k+1} / \Gamma_{k+2}$.  Now, suppose we are given $f, g \in \mcI_k$.
Then, we have $$f_*(g_*(x))x^{-1} = f_*(g_*(x)x^{-1})f_*(x)x^{-1} $$ $$= (f_*(t_g(x))t_g(x)^{-1}) t_g(x) f_*(x)x^{-1}
= t_f(t_g(x)) t_g(x) t_f(x)$$  Since $t_g(x) \in \Gamma_{k+1} / \Gamma_{k+2} \subseteq \krn{t_f}$, we find that
$f_*(g_*(x))x^{-1} = t_g(x) t_f(x)$
\end{proof}

\noindent \textbf{Remark:} In the above proof, we see that $\krn(t_f) \supset \Gamma_2 / \Gamma_{k+2}$, and so for $x
\in \Gamma_2 / \Gamma_{k+2}$ we have $$1 = t_f(x) = f_*(x)x^{-1} \Rightarrow f(x) = x$$  Thus $f$ acts trivially on
$\Gamma_2 / \Gamma_{k+2}$ and in particular on $\Gamma_{k+1} / \Gamma_{k+2}$.  Looking at the short exact sequence
\begin{equation} \label{eq:ses} 1 \to \Gamma_{k+1} / \Gamma_{k+2} \to \Gamma_1 / \Gamma_{k+2} \to
\Gamma_1 / \Gamma_{k+1} \to 1 \end{equation} one might think that $f$ must act trivially on $\Gamma_1 /\Gamma_{k+2}$
itself, but this is not the case.  Elements in $(\Gamma_1 / \Gamma_{k+2}) \setminus \, (\Gamma_2 / \Gamma_{k+2})$ may
be changed by elements in $\Gamma_{k+1} / \Gamma_{k+2}$ and this is precisely what $\tau_f$
measures. \\

In view of the remark, we see that $\tau_f$ retains the information of $f_* \in \Aut(\Gamma_1 / \Gamma_{k+2})$.
Furthermore, $\tau_f$ determines $f_*$ as an element of $\Aut(\Gamma_1 / \Gamma_{k+2})$ (assuming $f \in \mcI_k$). We
simply note that $f_*(x) = \tau_f(\overline{x}) x$ where $\overline{x}$ is the projection of $x$ to $H$.
Moreover, the following sequence is exact: (see Proposition 2.3 of \cite{M3})
\begin{equation} \label{eq:exactseq} 1 \to \Hom(H, \Gamma_{k+1} / \Gamma_{k+2}) \to
\Aut(\Gamma_1 /\Gamma_{k+2}) \to \Aut(\Gamma_1 / \Gamma_{k+1}) \end{equation}

Given $f \in \mcI_k$, one can similarly define a function $$\Gamma_m / \Gamma_{m+k+1} \to \Gamma_{m + k} /
\Gamma_{m+k+1}$$
\begin{equation} \label{eq:natural} x \mapsto f_*(x)x^{-1}
\end{equation}  As before, this induces a well-defined
homomorphism $\Gamma_m / \Gamma_{m+1} \to \Gamma_{m+k} / \Gamma_{m+k+1}$. (See Lemma 3.2 of \cite{M2}.)

Consider the free associative algebra $\mcA$ as defined in the previous section.  Suppose one has chosen
 $2g$ elements $\{x_1, ... , x_{2g}\} \subseteq \mcA$.  From general theory about the free associative algebra,
we know there is then a unique derivation $D: \mcA \to \mcA$ such that $D(a_i) = x_i$ where the $a_i$ are
generators of $\mcA$ (see \cite{R}). The following computation shows that $D(\mcL) \subseteq \mcL$ and that $D$ is a
derivation on $\mcL$:
$$ D[y, z] = D(yz - zy) = (Dy)z + yDz - (Dz)y - zDy = [Dy, z] + [y, Dz] $$

Thus, given $f \in \mcI_k$, there is a unique derivation $D_f$ of $\mcA$ which extends $\tau_f$.  It turns out
 that extending $\tau_f$ to all of $\mcL$ yields the same result regardless of whether one restricts
$D_f$ or uses (\ref{eq:natural}).  The following proposition follows more or less from Lemma 2.3 and Proposition
2.5 of \cite{M2}.

\begin{proposition}[Morita]
\label{proposition:tausagree} For all $m \geq 1$, the map defined by (\ref{eq:natural}) induces a homomorphism
$\Gamma_m / \Gamma_{m+1} \to \Gamma_{m+k} / \Gamma_{m+k+1}$ and is equal to the map $D_f\vert_{\mcL_m}$.
\end{proposition}

By abuse of notation, we will denote the extention to $\mcL$ by $\tau_f$.
The map $\tau$ has other nice algebraic properties.  They are collected in the following theorem.

\begin{theorem}[Morita]
Let $\tau$ be as defined above, a collection of homomorphisms $\mcI_k \to Der(\mcL)$, one for each $k$.  Then, the
following hold:
\begin{itemize}
\item[(a)] The map $\tau: \mcI_k \to Der(\mcL)$ is a homomorphism with kernel $\mcI_{k+1}(S)$.  Hence, it induces
a well-defined homomorphism $\mcI_k / \mcI_{k+1}(S) \to Der(\mcL)$
\item[(b)] The abelian group $$\bigoplus_{k=1}^\infty \mcI_k / \mcI_{k+1}(S) $$ has a Lie algebra structure induced by
$$\begin{array}{rcl} \mcI_m(S) \times \mcI_n(S) & \to & \mcI_{m+n}(S) \\ (f, g) & \mapsto & fgf^{-1}g^{-1} =: [f, g]
\end{array}$$

\item[(c)] $\tau$ induces a Lie algebra homomorphism $$ \bigoplus_{k=1}^\infty \mcI_k(S) / \mcI_{k+1}(S) \to Der(\mcL)$$
Furthermore, $\tau$ respects the conjugation action of $\Mod(S)$ on $\mcI_k$ and $Der(\mcL)$.
\end{itemize}
\end{theorem}
\begin{proof}[Proof Sketch:]  This proof sketch will consist mainly of citations.  For (a), recall that by
Proposition \ref{proposition:basictau}, $\tau_{f \circ g}|_H = \tau_f|_H + \tau_g|_H$.  Since the derivations
$\tau_{f \circ g}$ and $\tau_f + \tau_g$ agree on generators, they must agree on all of $\mcL$. One deduces the
kernel is $\mcI_{k+1}(S)$ from the exact sequence in (\ref{eq:exactseq}). Part (b) is Proposition 4.1 of \cite{M1}.
Also Proposition 4.7 of \cite{M1} shows (in slightly different notation) that $\tau_{[f,g]}|_H = (\tau_f \tau_g -
\tau_g \tau_f)|_H$. Since the two derivations $\tau_{[f,g]}$ and $\tau_f \tau_g - \tau_g \tau_f$ agree on $H$ and
since $H$ generates $\mcL$, we must have equality.  To show that the $\Mod(S)$ action is respected, we use the
definition of $\tau_f$ given by (\ref{eq:natural}).  Suppose $g \in \Mod(S)$.  In $\Gamma_m / \Gamma_{m+k+1}$
$$\begin{array}{rcl} \tau_{gfg^{-1}}(x) & = & g(f(g^{-1}(x))) x^{-1} = g(f(g(x)) g^{-1}(x^{-1}))\\ & = &g(f(g^{-1}(x))
(g^{-1}(x))^{-1}) = g(\tau_f(g^{-1}(x))) \end{array}$$
\end{proof}

\noindent \textbf{Remark:}  {\it A priori}, it may seem that, for $f \in \mcI_k$, we are using the entire action of $f_*$
on $\pi_1(S)$ since we use the action on $\Gamma_m / \Gamma_{m+k+1}$ for all $m$.  This would conflict with the
characterization given in the introduction that we only use the data of $f_* \in \Aut(\Gamma_1 /
\Gamma_{k+2})$.  However, since $\tau_f$ is a derivation on $\mcL$ which is generated by $H$, it is completely
determined by $\tau_f|_H$ which is itself determined by $f_* \in \Aut(\Gamma_1 / \Gamma_{k+2})$.

\section{The Contractions $\Phi_k$} \label{section:contractionsection}
 Our goal in
this section is to find a contraction ${\mathcal L}_{k+1} \to {\mathcal L}_1$ respecting the $\Sp$-action and thus the
$\Mod(S)$-action by the results of Section \ref{section:basicfacts}.
  We remark that we want to respect the action so that $\chi(\Psi_k(f))$ will depend only on the
conjugacy class of $f$ and because the argument in Section \ref{section:proofofcriterion} implicitly uses a change of
coordinates. The following theorem simplifies this problem. Below, $\Hom_{\Sp}$ will denote the set of homomorphisms which
respect the $\Sp$ action, and, for $X$ an $\Sp$-representation, $X_{\Sp}$ will indicate the space of vectors fixed by the
$\Sp$ action.  While I suspect the following may be known, I was not able to find it in the literature.

\begin{theorem}
\label{theorem:nondegen} If $f \in \Hom_{\Sp}({\mathcal L}_{k+1}, {\mathcal L}_1)$, then $\exists \phantom{|}n \in
\bold{Z}$ such that $nf$ is the restriction of an element $g \in \Hom_{\Sp}({\mathcal A}_{k+1}, {\mathcal A}_1)$,
where ${\mathcal A}_m$ is the summand $H^{\otimes m} \subset \mcA$.
\end{theorem}

\begin{proof}
The theorem will follow if we can find a bilinear pairing on each $\mcA_{k+1}$ which is nondegenerate on both
$\mcA_{k+1}$ and $\mcL_{k+1}$.
Let $\{a_1, b_1, ..., a_g, b_g\}$ be a symplectic basis of $H_1(S)$.  The $a_i$ and $b_i$
also serve as a free generating set of $\mcL$ as a Lie algebra and of $\mcA$ as an associative algebra.
We can easily define a pairing $\langle \; , \rangle$ which is nondegenerate on $\mcA_{k+1}$.  If $x = x_1 \otimes x_2 ...
\otimes x_{k+1}$ and $y = y_1 \otimes y_2 ... \otimes y_{k+1}$, then set
$$\langle x, y \rangle := \langle x_1, y_1 \rangle \langle x_2, y_2 \rangle ... \langle x_{k+1}, y_{k+1} \rangle$$
where $\langle x_i, y_i \rangle$ is the algebraic intersection pairing on $H$.

Now, let $\theta \in \Aut(\mcA)$ be the
algebra homomorphism defined by $\theta(a_i) = b_i$ and $\theta(b_i) = -a_i$.  In particular, if $w = x_1 \otimes x_2 \dots
 \otimes x_n$
then $\theta(w) = \theta(x_1) \otimes \theta(x_2) ... \otimes \theta(x_n)$. Let $Y_k$ be the canonical basis of $H^{\otimes
k}$ induced by the basis of $H$ (i.e. tensoring the $a$'s and $b$'s in every possible order). For two elements $y, y'
\in Y_k$, one easily sees that $\langle y, y' \rangle \neq 0$ if and only if $y' = \pm \theta(y)$. Then, for $P =
\sum_y c_y y$, we have $\langle P, \theta(P) \rangle > 0$, since all ``cross terms'' vanish and we are left with
$\sum_y c_y^2 \langle y, \theta(P) \rangle$.

We now wish to show that $\langle \, , \rangle$ is nondegenerate on
the embedded copy of ${\mathcal L}_{k+1}$, but this is almost immediate. We only need that $P \in
{\mathcal L}_{k+1}$ implies $\theta(P)\mcL_{k+1}$. Indeed, since ${\mathcal L}$ is the Lie
subalgebra of ${\mcA}$ generated by $\{a_1, b_1, ..., a_g, b_g\}$ and since $\theta$ preserves the Lie bracket and (up
to sign) permutes the generators $\{a_1, b_1, ..., a_g, b_g\}$, we see that $\theta(\mcL) = \mcL$.

Suppose $f \in \Hom_{\Sp}({\mathcal L}_{k+1}, {\mathcal L}_1) \cong ({\mathcal L}_{k+1}^* \otimes {\mathcal
L}_1)_{\Sp}$. Since $\mcL_{k+1}$ and $\mcL_{k+1}^*$ are finitely generated free $\bold{Z}$-modules, the pairing $\langle \; ,
\rangle$ gives an embedding ${\mathcal L}_{k+1} \hookrightarrow {\mathcal L}_{k+1}^*$ whose image has finite index.
Thus, there is some $n \in \bold{Z}$ such that $nf$ is in the image of $({\mathcal L}_{k+1} \otimes{\mathcal
L}_1)_{\Sp}$, but we have
$$({\mathcal L}_{k+1} \otimes{\mathcal L}_1)_{\Sp} \hookrightarrow (\mcA_{k+1} \otimes \mcA_1)_{\Sp}
\hookrightarrow (\mcA_{k+1}^* \otimes{\mcA}_1)_{\Sp} \cong \Hom_{\Sp}(\mcA_{k+1}, \mcA_1)$$ Thus, $nf$ is the
restriction of some $g \in \Hom_{\Sp}(\mcA_{k+1}, \mcA_1)$.
\end{proof}

Theorem \ref{theorem:nondegen} and its proof reduce our problem to finding tensors in $(\mcA_{k+1} \otimes
\mcA_1)_{\Sp} \cong (H^{\otimes k + 2})_{\Sp}$.  Thus, if $k = 2n$ is even, we obtain such a tensor by taking
the symplectic pairing $\omega_0 = \sum_i (a_i \otimes b_i - b_i \otimes a_i)$ and taking high tensor powers, i.e.
$\omega_0^{\otimes (n + 1)}$.  The element $\omega_0^{\otimes (n + 1)}$ represents the contraction
$$ x_1 \otimes x_2 ... \otimes x_{k+1} \mapsto \left( \prod_{j = 1}^{n} \langle x_{2j-1}, x_{2j}
\rangle \right) x_{k+1} $$ This contraction is what we denote by $\Phi_k$.

There is an obvious action of the permutation group $\mathfrak{S}_{2m}$ on $H^{\otimes 2m}$. Since $\Sp(2g, \Z)$ acts
diagonally on $H^{\otimes 2m}$, it is easy to see that for any $\sigma \in \mathfrak{S}_{2m}$, we have $\eta \in
(H^{\otimes 2m})_{\Sp}$ if and only if $\sigma(\eta) \in (H^{\otimes 2m})_{\Sp}$.  Thus, all the vectors $\sigma
(\omega_0^{\otimes 2m})$ are $\Sp$-invariant as well.  For every $\sigma \in \mathfrak{S}_{2m}$, there is a
corresponding $\sigma'$ so that $\sigma(\omega_0^{2m})$ corresponds to the contraction
$$ x_1 \otimes x_2 ... \otimes x_{m-1} \mapsto \left( \prod_{j = 1}^{n} \langle x_{\sigma'(2j-1)}, x_{\sigma'(2j)}
\rangle \right)x_{\sigma'(m-1)}$$  Furthermore, it is a classical result of Weyl (see, e.g., Section 4.1 of \cite{M4}) that
$\displaystyle \{\sigma(\omega_0^{\otimes 2m})\}_{\sigma \in \mathfrak{S}_{2m}}$ is a generating set for $((H \otimes
\Q)^{\otimes 2m})_{\Sp(2g, \Q)}$.

\section{Proof of theorem \ref{theorem:maintheorem}}
\label{section:proofofcriterion}

\noindent Recall from above that for
each $k \geq 1$ we defined a map $$\begin{array}{rcl} \Psi_k: \mcI_k & \to & \End(H) \\
f & \mapsto & \left\{ \begin{array}{ll} \Phi_k \circ (\tau_f|_H)  & \textrm{$k$ even} \\
                                         \Phi_{2k} \circ (\tau_f^2|_H)  & \textrm{$k$ odd}  \end{array} \right. \end{array}
                                         $$

We remark that the following proof of the main theorem remains valid if we replace $\Phi_{k}$ with any of the
contractions induced by a $\sigma(\omega_0^{k+2})$ described in Section \ref{section:contractionsection}. In the
following, all factorization and irreducibility is with respect to $\Z[x]$.

\begin{proof}[Proof of Theorem \ref{theorem:maintheorem}]
Let $f \in \mcI_k$.  Recall that the Nielsen--Thurston classification and torsion-freeness of $\mcI_1 \supseteq
\mcI_k$ imply that $f$ is pseudo-Anosov if and only if $f$ is irreducible.  It is well-known that $\mcI_1$
is {\em pure}, meaning that if an isotopy class of 1-submanifold is fixed, then each component of the 1-submanifold
is fixed (see Theorem 1.2 of \cite{I}).  Thus, the proof of Theorem \ref{theorem:maintheorem} reduces to proving the
 following two claims. \\

\textbf{Claim 1:} Suppose $f$ fixes an essential separating curve. Then, the characteristic polynomial of $\Psi_k(f)$
factors into two (nontrivial) even degree polynomials in $\Z[x]$. \\

\textbf{Claim 2:} Suppose $f$ fixes a nonseparating curve.  Then, $\Psi_k(f)$ has an eigenvector over $\Z$. \\

 \noindent

Before we begin the proofs of Claims 1 and 2, we state a theorem that will be used for both.  (This is Theorem 2.5 in
\cite{R})
\begin{theorem}[Shirshov, Witt]
\label{theorem:subalgfreealg} If ${\mathcal L}'$ is a subalgebra of a free Lie algebra ${\mathcal L}$ over a field,
then ${\mathcal L}'$ is a free Lie algebra .
\end{theorem}

\noindent \textbf{Proof of Claim 1:}  Let $\gamma$ be the (oriented) separating curve such that $f(\gamma) = \gamma$.
Cutting along $\gamma$ separates $S$ into a $\Sigma_{g_1, 1} =: S_1$ and a $\Sigma_{g_2, 2} =: S_2$ where $g_1 + g_2
= g$. Let $C$ (resp. $D$) be the image of $H_1(S_1, \Z)$ (resp. $H_1(S_2, \Z)$) in $H$.  Since, $f(S_i) = S_i$ (up to
isotopy), one might hope that either $\Psi_k(f)(C) \subseteq C$ or $\Psi_k(f)(D) \subseteq D$. We will show that this
actually holds for $D$.

We begin by defining a submodule of $\mcL$:
$$ M := \bigoplus_m ( \Lambda \cap \Gamma_m / \Lambda \cap \Gamma_{m+1} ) $$ where $\Lambda := \pi_1(S_2)$.
Note that $M \cap \mcL_1 = D$.  {\em Step 1} is to show that $\tau_f(M) \subseteq M$. {\em Step 2} is to show that
$M$ is a free Lie subalgebra and give generators of $M$
as a Lie algebra. {\em Step 3} is to show,
using the generators, that for any $x \in M$ we have $\Phi_n(x) \in D$. Then, it is clear from the definition of
$\Psi_k$ that for $d \in D$, we have $\Psi_k(f)(d) \in D$.  Since $D$ is an even rank subspace, that will complete
the proof.

First, we need to set up some notation.  Let $p_1 \in \partial S_1$ (resp. $p_2 \in \partial S_2$) be the basepoint
of $S_1$ (resp $S_2$ and $S$). Let
$\alpha$ be a path from $p_2$ to $p_1$, and let $\tilde{\gamma} = \alpha \gamma \alpha^{-1} \in \pi_1(S_2)$.  Let
$\iota$ (resp. $\hat{\iota}$) denote geometric (resp. algebraic) intersection number of unbased homotopy classes of
closed curves. Choose $\{c_i'\}_{i=1}^{2g_1}
\in \pi_1(S_1, p_1)$ and $\{d_i\}_{i=1}^{2(g_2)} \in \pi_1(S_1, p_2)$ with the following properties (see Figure 1): \\
\begin{itemize}
\item[(a)] The set $\{c_i'\}_{i=1}^{2g_1}$ (resp. $\{\tilde{\gamma}\} \cup \{d_i\}_{i=1}^{2(g_2)}$) generates $\pi_1(S_1,
p_1)$ (resp. $\pi_1(S_2, p_2)$). \\
\item[(b)] For all $m, n$, we have $\iota(c_m', d_n) = \hat{\iota}(c_m', d_n) = 0$.  Furthermore,
$$ \begin{array}{rcl}
 \iota(c_m', c_n') & = & \left\{ \begin{array}{ll} 1 & \textrm{if } m = n + g_1 \textrm{ or } m = n - g_1 \\
                                            0 & otherwise \end{array} \right. \\
 \iota(d_m, d_n) & = & \left\{ \begin{array}{ll} 1 & \textrm{if } m = n + g_2 \textrm{ or } m = n - g_2 \\
                                            0 & otherwise \end{array} \right.
\end{array}$$
and for $1 \leq i \leq g_1$ (resp. $1 \leq i \leq g_2$), we have $\hat{\iota}(c_i', c_{i+g_1}') = 1$ (resp.
$\hat{\iota}(d_i, d_{i+(g_2)}) = 1$ \\
\item[(c)] As an element of $\pi_1(S_1, p_1)$, we have $\gamma = \prod_{i=1}^{g_1} [c_i', c_{i+g_1}']$. \\
\end{itemize}

In particular, the union $\{c_i'\}_{i=1}^{2g_1} \cup \{d_i\}_{i=1}^{2g_2}$ gives a symplectic basis in $H$. Now, let
$c_i := \alpha c_i' \alpha^{-1}$. We have $\tilde{\gamma} = \prod_{i=1}^{g_1} [c_i, c_{i+{g_1}}]$ and $\pi_1(S, p_2) =
\langle \{ c_i \}, \{d_i\} \rangle $. Furthermore, denote the inclusion map of $S_2$ by $j: S_2 \hookrightarrow S$.
In the following, we will frequently view $d_i \in \mcL_1$ and $\tilde{\gamma} \in \mcL_2$. \\

 \begin{figure}[htbp]
 \begin{center}
 \psfrag{cg1+1'}{$c_{g_1+1}'$}
 \psfrag{c1'}{$c_1'$}
 \psfrag{c2g1'}{$c_{2g_1}'$}
 \psfrag{cg1'}{$c_{g_1}'$}
 \psfrag{p1}{$p_1$}
 \psfrag{d1}{$d_1$}
 \psfrag{dg2+1}{$d_{g_2+1}$}
 \psfrag{d2g2}{$d_{2g_2}$}
 \psfrag{dg2}{$d_{g_2}$}
 \psfrag{p2}{$p_2$}
 \psfrag{S1}{$S_2$}
 \psfrag{S2}{$S_2$}
 \psfrag{gamma}{$\gamma$}
 \psfrag{...}{$\dots$}
 \includegraphics{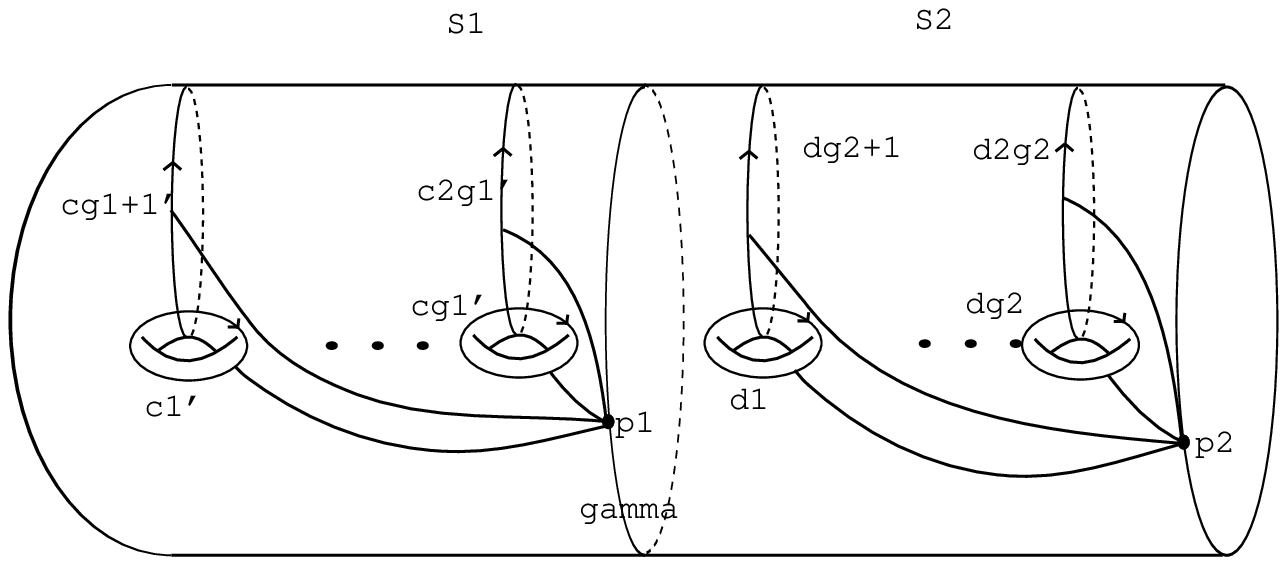}

 \caption{}
 \end{center}
 \end{figure}

\noindent {\em Step 1:} First note that since $S$ and $S_2$ share a base point, $\pi_1(S_2)$ gives a well-defined
subgroup of $\pi_1(S) = \Gamma$ which is invariant under $f_*$. We remark that a similar statement is not true for
$S_1$. Indeed, to embed $\pi_1(S_1)$ in $\pi_1(S)$ requires that we choose a path connecting base points (e.g.
$\alpha$); even after choosing a representative homeomorphism of $f$ which fixes $\gamma$ pointwise, this path is not
necessarily preserved (up to homotopy rel endpoints).

 Recall that one way of defining $\tau_f$ is to induce it from the map
$$\begin{array}{rcl} \Gamma_m & \to & \Gamma_{m+k} \\ x & \mapsto & f_*(x)x^{-1} \end{array}$$
Since $f_*(\Lambda) = \Lambda$, it is easy to see that
$$ M = \bigoplus_m ( \Lambda \cap \Gamma_m / \Lambda \cap \Gamma_{m+1} ) $$ is a $\tau_f$-invariant submodule of $\mcL$. \\

\noindent {\em Step 2:} We wish to show $M$ is a Lie subalgebra and find its generators.  We will do this showing
that $M$ is the Lie algebra homomorphic image of a Lie algebra $N$ whose generators are easily found.

We first define a filtration of $\Lambda$ which is a slight alteration of the lower central series.   We let
$$\begin{array}{lll} \Lambda_1 := \pi_1(S_2) & \Lambda_2 := \langle [\Lambda_1, \Lambda_1], \tilde{\gamma} \rangle
& \Lambda_m := \langle [\Lambda_{m-n}, \Lambda_n] \rangle_{n=1}^{\lfloor{\frac{m}{2}}\rfloor} \;\;\; \textrm{for} \;
m
 \geq 3 \end{array}$$
By Theorem \ref{theorem:liealg}, $$ N := \bigoplus_n \Lambda_n / \Lambda_{n+1}$$ is a graded $\Z$-Lie algebra under
the commutation bracket. Since $j_*(\Lambda_n) \subseteq \Gamma_n$, there is an induced Lie algebra homomorphism $N \to
\mcL$. It is easy to check that, as a Lie algebra, $N$ is generated by $\{d_i\}_{i=1}^{2g_2} \cup
\{\tilde{\gamma}\}$ and so its image $M' := j_*(N)$ in $\mcL$ is also generated by $\{d_i\}_{i=1}^{2g_2} \cup
\{\tilde{\gamma}\}$ (viewed in $\mcL$).

\begin{proposition}
\label{proposition:injective} $N$ maps isomorphically onto $M'$
\end{proposition}
\begin{proof}[Proof of Proposition \ref{proposition:injective}]

We wish to use Theorem \ref{theorem:subalgfreealg}, but $\mcL$ is not an algebra over a field.  As $\Q$ is a flat
$\Z$-module, we have $M' \otimes \Q \hookrightarrow \mcL \otimes \Q$, and so $M_{\Q}' := M' \otimes \Q$ is a free Lie
algebra generated by $\{d_i\}_{i=1}^{2(g_2)} \cup \{\tilde{\gamma}\}$, but it is not a priori clear that these
generators are {\em free}. In the proof of Theorem \ref{theorem:subalgfreealg} in \cite{R}, a recipe is given for finding
free generators of a subalgebra, which we describe now.

For any subset $X \subseteq \mcL \otimes \Q$, let $\langle X \rangle$ denote the Lie subalgebra of $\mcL \otimes \Q$
generated by $X$.  Let $$E_n = M_{\Q}' \cap \left( \bigoplus_{i=1}^n \mcL_i \otimes \Q \right) $$ and let $E_n' = E_n
\cap \langle E_{n-1} \rangle $.  If we let $X_n := $ a set of generators (as a $\Q$ vector space) for $E_n$ mod
$E_n'$, then $X = \bigcup_n X_n$ is a free generating set of $M_{\Q}'$.

We now show the afore-mentioned generators of $M_{\Q}'$ to be free.  Clearly, we can set $X_1 := \{d_i\}_{i=1}^{2g_2}$.
The only question is whether $\tilde{\gamma}$ is in the Lie algebra generated by $X_1$. Recall that $\tilde{\gamma} =
\prod_{i=1}^{g_1} [c_i, c_{i+g_1}]$ and so in $\mcL_2$, we have $\tilde{\gamma} = \sum_{i=1}^{g_1} [c_i, c_{i+g_1}]$.
As elements of $H$, the $c_i$ and $d_i$ freely generate $\mcL \otimes \Q$, so $\tilde{\gamma} \notin \langle X_1
\rangle$. Thus, we can set $X_2 = \{\tilde{\gamma}\}$, and so $\{d_i\}_{i=1}^{2g_2} \cup \{\tilde{\gamma}\}$ freely
generates $M_{\Q}'$.  But then clearly it {\em freely} generates $M'$.

Now, we can define an inverse Lie homomorphism $M' \to N$ by sending generators to generators, and so $N \to \mcL$ is
injective.
\end{proof}
By the proposition, we have $\Lambda_n \setminus \Lambda_{n+1} \hookrightarrow \Gamma_n \setminus \Gamma_{n+1}$, but
this implies that in fact $\Lambda_n = \Lambda \cap \Gamma_n$.  Thus, $M = M'$. \\

{\em Step 3:} Recall that $C = $ image of $H_1(S_1)$ and $D = $ image of $H_1(S_2)$ in $H$; i.e. $C = \langle
\{c_i\}_{i=1}^{2g_1} \rangle$ and $D = \langle \{d_i\}_{i=1}^{2(g_2)} \rangle$ . Suppose $x \in D$. Then, by Steps
1 and 2,
$$y := \left\{ \begin{array}{ll} \tau_f(x)  & \textrm{$k$ even} \\
                                         \tau_f^2(x)  & \textrm{$k$ odd}  \end{array} \right. $$
is an element of $M$.  We can write $\tilde{\gamma}$ in $\mcA$ as $\sum_{i=1}^{2g_1} (c_i \otimes c_{i+g_1} -
c_{i+g_1} \otimes c_i)$. Thus, $M$ is contained in the subring generated by
$$\left\{\sum_{i=1}^{2g_1} (c_i \otimes c_{i+g_1} - c_{i+g_1} \otimes c_i)\right\} \cup \{d_i\}_{i=1}^{2(g_2)}$$
 Consequently, we can write $y = \sum_m y_{m, 1} \otimes ... \otimes y_{m, n}$ where an
{\em even} number of the elements of $\{y_{m, 1}, ..., y_{m, n}\}$ are in $C$ and the rest are in $D$. Since
$\hat{\iota}(c_i, d_j) = 0$ for all $i, j$, we have $\Phi_{n-1}(y_{m, 1} \otimes ... \otimes y_{m, n}) \neq 0$ only
if $y_{m, n} \in D$.  Thus, $\Psi_k(f)(D) \subseteq D$, and we are done with Claim 1.\\

\noindent \textbf{Proof of Claim 2:} Let $\alpha$ be the nonseparating curve which is fixed by $f \in \mcI_k$.  Let
$\hat{S}$ be the surface obtained by cutting along $\alpha$, and $j: \hat{S} \hookrightarrow S$ the canonical
immersion. Similar to the proof of Claim 1, we will show that $\Psi_k(f)(C) \subseteq C$ where $C := $ image of
$H_1(\hat{S}, \Z)$ in $H$.  Analagous to the above, we let $$M := \bigoplus_n \hat{\Gamma} \cap \Gamma_n /
\hat{\Gamma} \cap \Gamma_{n+1}$$ where $\hat{\Gamma} = \pi_1(\hat{S})$.  We go through the same 3 steps as in
the proof of Claim 1:
\begin{itemize}
\item {\em Step 1:} Show that $\tau_f(M) \subseteq M$.
\item {\em Step 2:} Show that $M$ is a Lie subalgebra of $\mcL$ and find generators.
\item {\em Step 3:} Show that $\Phi_k(M \cap \mcL_{k+1}) \subseteq C$.
\end{itemize}

Let us first set up some notation.  Let $\alpha_1$ and $\alpha_2$ be the boundary curves of $\hat{S}$ such that
$j_*(\alpha_1) = j_*(\alpha_2) = \alpha$. Choose based representatives $a, a_1$ and $a_2$ of $\alpha, \alpha_1$ and
$\alpha_2$ respectively as in Figure 2; in particular, $j_*(a_1) = a$. Also, let $b$ be as depicted in
Figure 2.  Extend $\{a, b\}$ to a
``standard'' generating set $\{a, b\} \cup \{c_i\}_{i=1}^{2(g-1)}$; i.e. the following hold: \\

\begin{itemize}
\item[(a)]  The set $\{a, b\} \cup \{c_i\}_{i=1}^{2(g-1)}$ gives a symplectic basis in homology.
\item[(b)]  $\iota(a,b) = \hat{\iota}(a, b) = 1$.
\item[(c)]  All $c_i$ can be homotoped to lie entirely inside the interior of $\hat{S}$.
\end{itemize}
Letting $a_1$ and $a_2$ be as in Figure 2, one can easily check that  $j_*(a_1 a_2^{-1}) = [a, b^{-1}]$. \\

\begin{figure}[htbp]
\begin{center}
\psfrag{a}{$a$} \psfrag{b}{$b$} \psfrag{alpha}{$\alpha$} \psfrag{S}{$S$} \psfrag{Shat}{$\hat{S}$}
\psfrag{alpha1}{$\alpha_1$} \psfrag{alpha2}{$\alpha_2$} \psfrag{a1}{$a_1$} \psfrag{a2}{$a_2$}

\includegraphics{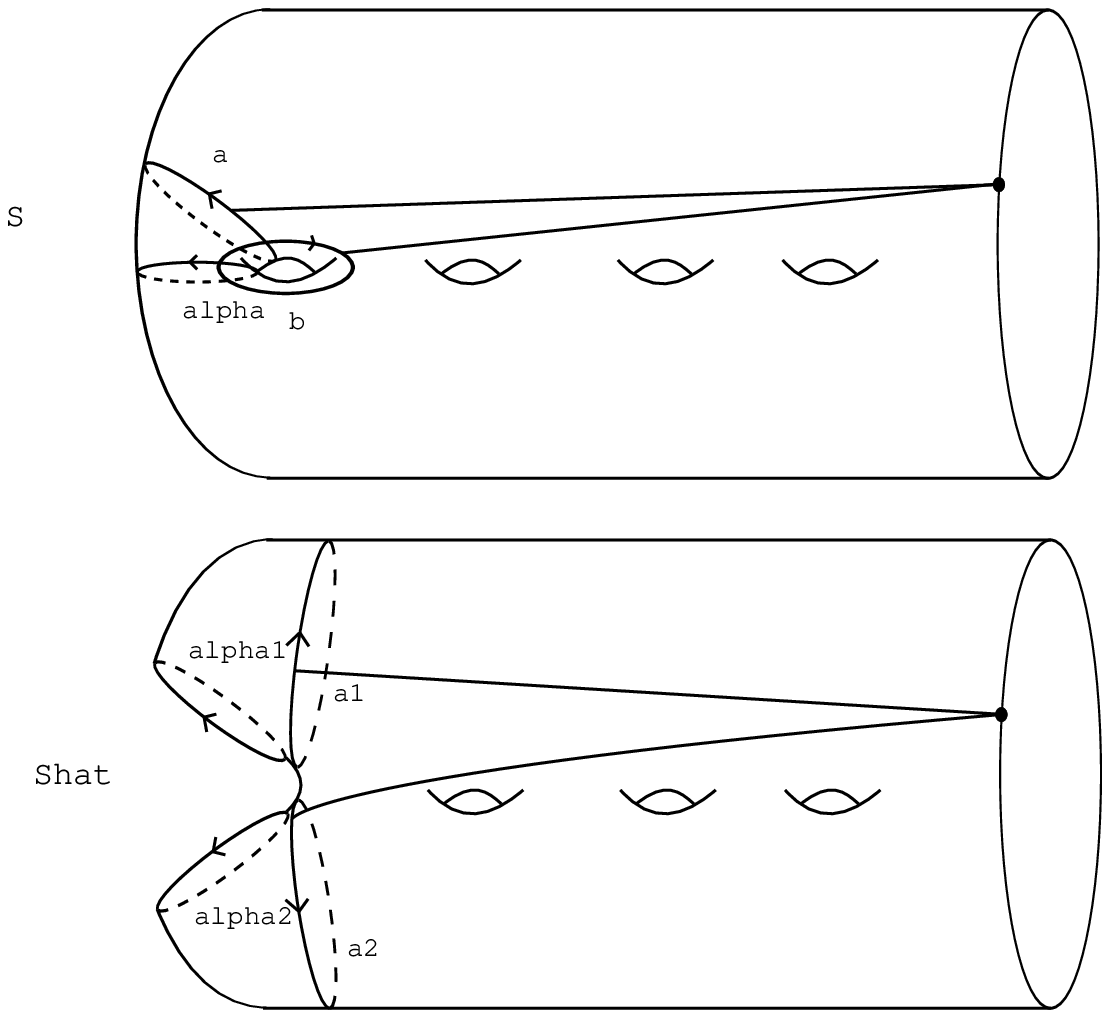}
\caption{}
\end{center}
\end{figure}

\noindent {\em Step 1:} Choosing the same basepoint for $\hat{S}$ and $S$, we have that
$j_*: \pi_1(\hat{S}) \to \pi_1(S)$ is injective and $\hat{\Gamma} = \pi_1(\hat{S})$ is invariant under
$f_*$.  Thus, we have
$$M := \bigoplus_n \hat{\Gamma} \cap \Gamma_n / \hat{\Gamma} \cap \Gamma_{n+1}$$
is a $\tau_f$-invariant submodule of $\mcL$.  It is also easy to see $M \cap \mcL_1 = C$. \\

\noindent {\em Step 2:}
 Just as in the proof of Claim 1, we choose a filtration of
$\pi_1(\hat{S})$ which is a slight alteration of the lower central series:
$$\begin{array}{rcl}\hat{\Gamma}_1 & = & \pi_1{\hat{S}} \\ \hat{\Gamma}_2 & = & \langle [\hat{\Gamma}_1, \hat{\Gamma}_1], a_1
a_2^{-1} \rangle \\  \hat{\Gamma}_n & = & \langle [\hat{\Gamma}_{n-k}, \hat{\Gamma}_k]
\rangle_{k=1}^{\lfloor{\frac{m}{2}}\rfloor} \;\; n \geq 3 \end{array}$$   By Theorem \ref{theorem:liealg}, we get a
corresponding graded $\Z$-Lie algebra which we denote by $\hat{M}$.  Again, since $j_*(\hat{\Gamma}_n) \subseteq
\Gamma_n$, we get an induced Lie algebra homomorphism $\hat{M} \to \mcL$.  Note that $\hat{M}$ is generated by
$\{a_1\} \cup \{c_i\}_{i=1}^{2(g-1)} \in \hat{M}_1$ and $a_1 a_2^{-1} \in \hat{M}_2.$  Since $a_1 a_2^{-1} \mapsto
[a, b^{-1}]$, we have that $\{a, [a, b^{-1}]\} \cup \{c_i\}_{i=1}^{2(g-1)}$ generates $j_*(\hat{M})$.
\begin{proposition}
\label{proposition:injective2} The Lie algebra $\hat{M}$ maps isomorphically onto $j_*(\hat{M})$.
\end{proposition}
\begin{proof}[Proof of Proposition \ref{proposition:injective2}]
Since the set $\{a, b\} \cup \{c_i\}_{i=1}^{2(g-1)}$ is a free generating set of $\mcL$, we have $[a, b^{-1}] \notin
\langle a, \{c_i\}_{i=1}^{2(g-1)} \rangle$.  Thus, by reasoning similar to that in the separating case, $\{a, [a,
b^{-1}]\} \cup \{c_i\}_{i=1}^{2(g-1)}$ is a {\em free} generating set of $j_*(\hat{M})$. We obtain an inverse Lie
algebra map $j_*(\hat{M}) \to \hat{M}$ induced by $$\begin{array}{rcl} a & \mapsto & a_1 \\ \phantom{a} [a, b^{-1}] &
\mapsto & a_1 a_2^{-1} \\ c_i & \mapsto & c_i \end{array}$$
\end{proof}
Since $\hat{M}$ injects into $\mcL$, we have
$$\hat{\Gamma}_m \setminus \hat{\Gamma}_{m+1} \hookrightarrow \Gamma_m \setminus \Gamma_{m+1}$$ and so
$\hat{\Gamma}_m = \hat{\Gamma} \cap \Gamma_m$. Thus, $j_*(\hat{M}) = M$.\\

{\em Step 3:} \noindent Now, let $x \in C:= \langle a, \{c_i\}_{i=1}^{2(g-1)} \rangle \subseteq H$.  Then
$$y := \left\{ \begin{array}{ll} \tau_f(x)  & \textrm{$k$ even} \\
                                         \tau_f^2(x)  & \textrm{$k$ odd}  \end{array} \right. $$
is an element of $M$.  As an element of $\mcA$, we may write $y = \sum_m y_{m, 1} \otimes ... \otimes y_{m, n}$ where
each $y_{m, r}$ is a multiple of one of $a, b, c_i$.  Since $[a, b] = a \otimes b - b \otimes a$ as an element of
$\mcA$ and $y \in
M$, there are at least as many $a$ terms as $b$ terms in $y_{m, 1}, ...,  y_{m, n}$.  Since $b$ pairs nontrivially
only with $a$ in the set $\{a, b\} \cup \{c_i\}_{i=1}^{2(g-1)}$, we have $\Phi_{n-1}(y_{m, 1} \otimes ... \otimes
y_{m, n}) \neq 0$ only if $y_{m, n} \neq$ a multiple of $b$, in which case $\Phi_{n-1}((y_{m, 1} \otimes ... \otimes
y_{m, n})\in C$.  Thus, $\Psi_k(f)(C) \subseteq C$, and since $C$ has rank
$2g - 1$, the characteristic polynomial of $\Psi_k(f)$ factors into a product of a degree $1$ and degree $2g - 1$ polynomial.
\end{proof}

\section{Theorem \ref{theorem:maintheorem} vs. the Thurston--Penner Criteria}
\label{section:example} \noindent In this section we will compare the criterion of Theorem \ref{theorem:maintheorem}
to the Thurston--Penner criteria.  Since the Thurston--Penner criteria are topological and Theorem
\ref{theorem:maintheorem} is algebraic, one might expect that there is essentially no relation between the two.  We
will show this to be true in the following sense. There exist examples satisfying the Thurston or Penner criteria but
not the hypothesis of Theorem \ref{theorem:maintheorem} and examples satisfying both.  As of the writing of this
paper, it has not been proven that there are examples of pseudo-Anosovs which do not satisfy the Thurston--Penner
criteria. However, we will give an example satisfying the hypothesis of Theorem \ref{theorem:maintheorem} to which
the Thurson--Penner criteria do not seem to apply directly.

Since we will be dealing with Dehn twists about separating curves, we first describe $\Psi_2(T_{\gamma})$ where
$\gamma$ is a ``standard'' separating curve and $T_\gamma$ is the Dehn twist about $\gamma$.  First let us set up a
symplectic basis.  Let $\{\alpha_i, \beta_i\}$ be the curves as depicted in Figure 3 with $a_i = [\alpha_i]$ and $b_i
= [\beta_i]$ their homology classes. Our ordered basis of $H$ throughout this section will be $\{a_1, b_1, a_2, b_2,
\dots, a_g, b_g\}$. By ``standard'' separating curve, we will mean one of the $\gamma_i$ as depicted in Figure 3.

 \begin{figure}[htbp]
 \begin{center}
 \psfrag{beta1}{$\beta_1$}
 \psfrag{beta2}{$\beta_2$}
 \psfrag{betai}{$\beta_i$}
 \psfrag{betai+1}{$\beta_{i+1}$}
 \psfrag{betag}{$\beta_g$}
 \psfrag{alpha1}{$\alpha_1$}
 \psfrag{alpha2}{$\alpha_2$}
 \psfrag{alphai}{$\alpha_i$}
 \psfrag{alphai+1}{$\alpha_{i+1}$}
 \psfrag{alphag}{$\alpha_g$}
 \psfrag{gamma1}{$\gamma_1$}
 \psfrag{gammai}{$\gamma_i$}
 \psfrag{gammag-1}{$\gamma_{g-1}$}
 \psfrag{...}{$\dots$}

 \includegraphics{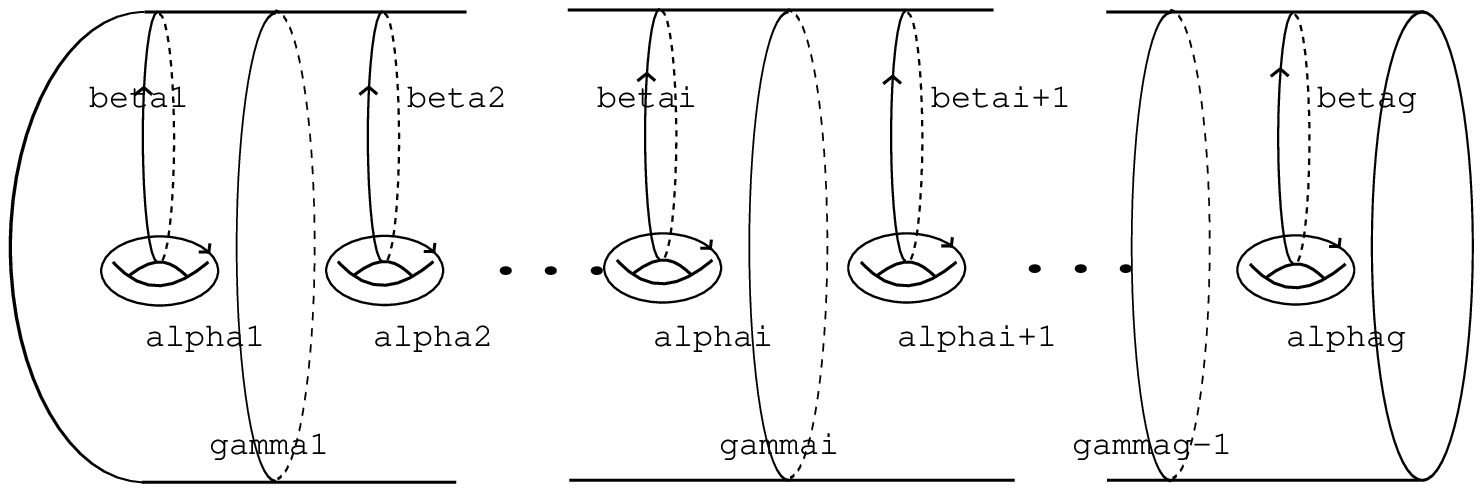} \caption{}
 \end{center}
 \end{figure}

\begin{lemma} \label{lemma:psisepcurve}
With $\{a_i, b_i\}$ and $\{\gamma_i\}$ as above, the element $\Psi_2(T_{\gamma_i}) \in \End(H)$ is the map defined
by:
$$\begin{array}{rcl} a_j & \mapsto & \left\{ \begin{array}{ll} (2i + 1) a_j & j \leq i \\ 0 & j > i \end{array} \right. \\
                    b_j & \mapsto & \left\{ \begin{array}{ll} (2i + 1) b_j & j \leq i \\ 0 & j > i \end{array} \right.
                    \end{array} $$
\end{lemma}

\noindent \textbf{Remark:} Note that with the given indexing, $i$ is the genus of $\gamma_i$.

\begin{proof}
We can lift $a_i, b_i, \gamma_i$ to $\tilde{a_i}, \tilde{b_i}, \tilde{\gamma_i} \in \pi_1(S)$ by connecting
$\alpha_i$, $\beta_i$ and $\gamma_i$ to the basepoint via paths.  Furthermore, we can do it in such a way that
 $\tilde{\gamma_i} = \prod_{j=1}^i [\tilde{a_j}, \tilde{b_j}]$ in $\pi_1(S)$ and
$$\begin{array}{rcl} T_{\gamma_i}(\tilde{a_j}) & = & \left\{ \begin{array}{ll} \tilde{\gamma_i} \tilde{a_j}
\tilde{\gamma_i}^{-1} & j \leq i \\ \tilde{a_j} & j > i \end{array} \right. \\
                   T_{\gamma_i}(\tilde{b_j}) & = & \left\{ \begin{array}{ll} \tilde{\gamma_i} \tilde{b_j}
\tilde{\gamma_i}^{-1} & j \leq i \\ \tilde{b_j} & j > i \end{array} \right.    \end{array} $$

Thus, for $j \leq i$ and $f_i = T_{\gamma_i}$, we compute $f_i(\tilde{a_j}) \tilde{a_j}^{-1} = [\tilde{\gamma_i},
\tilde{a_j}]$ and
$$\tau_{f_i} (a_j) = [\sum_{k=1}^i [a_k, b_k], a_j] = \sum_{k=1}^i
((a_k \otimes b_k - b_k \otimes a_k) \otimes
 a_j - a_j \otimes (a_k \otimes b_k - b_k \otimes a_k))$$

For $j > i$, we easily see that $\tau_{f_i}(a_j) = 0$.  Recall that $\Phi_2(c_1 \otimes c_2 \otimes c_3) =
\hat{\iota}(c_1, c_2) c_3$. We then compute for $j \leq i$ that $\Psi_2(f_i) = \Phi_2 (\tau_{f_i}(a_j)) = (2i + 1)
a_j$.  Clearly, $\Phi_2(\tau_{f_i}(a_j)) = 0$ for $j > i$.  The computation for $b_j$ is the same but with the
the roles of $a$ and $b$ switched.
\end{proof}

Now let us consider $T_{\gamma}$ where $\gamma$ is an arbitrary separating curve not homotopic to the boundary.
Recall that $\Psi_k$ is $\Mod(S)$-equivariant (This follows from the $\Mod(S)$-equivariance of $\Phi_k$ and $\tau$).
The $\Mod(S)$ action on $\End(H)$ is as follows. If $\varphi \in \Mod(S)$
 and $h \in \End(H)$, then $$\varphi \cdot h =
[\varphi] h [\varphi]^{-1}$$  where $[\varphi]$ denotes the projection of $\Phi$ to $\Sp(2g, \Z)$.
 Thus, for $f \in \mcI_2$ and $\varphi \in \Mod(S)$, we find that $\Psi_k(\varphi f \varphi^{-1}) = [\varphi] \Psi_k(f)
[\varphi]^{-1}$. Recall that if for a fixed $g'$, two separating curves $\eta_1$ and $\eta_2$ both cut $S$ into a
$\Sigma_{g', 1}$ and a $\Sigma_{g-g', 2}$, then there is some $\varphi \in \Mod(S)$ such that $\varphi(\eta_1) = \eta_2$.  Thus,
$\Psi_2(T_{\gamma})$ is of the form $\varphi \Psi_2(T_{\gamma_i}) \varphi^{-1}$ for some $i$ and some $\varphi \in \Sp(2g, \Z)$.
Similarly, if $A$ is a multicurve of separating curves and $T_A$ the multicurve twist, then $$ \Psi_2(T_A) = \varphi
\Psi_2(\prod_{k=1}^m T_{\gamma_{i_k}}) \varphi^{-1}$$ for some $\varphi \in \Sp(2g, \Z)$ and some subset $\{\gamma_{i_k}\}$ of
$\{\gamma_i\}$.

For the reader's convenience, we recall a few definitions and state a corollary to both the Thurston and Penner
criteria.  A {\em pants decomposition} is a maximal set of pairwise nonisotopic simple closed curves which are
pairwise disjoint not null-homotopic.  For an $S_{g,b}$, a pants decomposition consists of $3g - 3 + 2b$ curves.
Recall that a simple closed curve $\gamma$ is {\em essential} if it is neither homotopically trivial
nor homotopic to a
boundary component.  We say that two curves $\eta$ and $\nu$ {\em fill} a surface $S$ if, for any essential simple
closed curve $\gamma$, the curve $\gamma$ either intersects $\eta$ or $\nu$ nontrivially.  We define the notion of
filling for two multicurves similarly.
\begin{corollary}[Thurston, Penner] \label{corollary:thustonpenner}
If two multicurves $A$ and $B$ fill a surface, then the product of multicurve twists $T_A T_B^{-1}$ is pseudo-Anosov.
\end{corollary}

\subsection{Negative Results for Theorem \ref{theorem:maintheorem}}
In this section we show that there is a pseudo-Anosov in $\mcI_2(S_{g,1})$ for each $g \geq 2$ which satisfies the
Thurston--Penner criteria but not the hypothesis of Theorem \ref{theorem:maintheorem}.  Let $T_\gamma$ denote the
twist about a simple closed curve $\gamma$.
\begin{theorem}
For each $g \geq 2$, there exists two simple closed curves $\gamma_{g, 1}$ and $\gamma_{g, 2}$ filling $S = S_{g,1}$
such that $f_g := T_{\gamma_{g,1}} T_{\gamma_{g,2}}^{-1}$ does not satisfy the hypothesis of Theorem
\ref{theorem:maintheorem}.  However, by the Thurston--Penner criteria, we know $f_g$ is pseudo-Anosov.
\end{theorem}
\begin{proof}
We break the proof into two cases.  For $g = 2$, we will explicitly compute $\Psi_2(f_2)$.  For $g \geq 3$, we will
show that there is an $f_g'$ such that $f_g'$ is reducible and $\Psi_2(f_g') = \Psi_2(f_g)$.  Of course, then it is
impossible for $\Psi(f_g)$ to satisfy the hypothesis of Theorem \ref{theorem:maintheorem} since $\Psi(f_g')$ does
not.

We also need a consequence of Lemma 2 of Expose 13 of \cite{FLP} to construct the $f_g$.  For the reader's
convenience, we state the consequence.
\begin{lemma} \label{lemma:flp}
Let $S$ be a surface.  Let $\gamma$ be a simple closed curve on $S$ and $P = \{ \alpha_1, \dots, \alpha_{3g-3} \}$ a
pants decomposition of $S$ such that $\iota(\gamma, \alpha_i) \neq 0$ for all $\alpha_i$ that are not boundary
components. Then, the curves $\gamma$ and $T_P(\gamma)$ fill the surface.
nontrivially.
\end{lemma}

\noindent \textbf{Case $g = 2$:}  Let $\gamma_{2,1}$ and the $\eta_i$ be as in Figure 4.  Since the $\eta_i$ are
disjoint and $\{\eta_i\}$ is a 4 element set, $P = \{\eta_i\}$ is a pants decomposition.  By Lemma \ref{lemma:flp},
we know that $\gamma_{2,1}$ and $ \gamma_{2,2} := T_P(\gamma)$ fill $S$.

 \begin{figure}[htbp]
 \begin{center}
 \psfrag{gamma}{$\gamma$}
 \psfrag{eta1}{$\eta_1$}
 \psfrag{eta2}{$\eta_2$}
 \psfrag{eta3}{$\eta_3$}

 \includegraphics{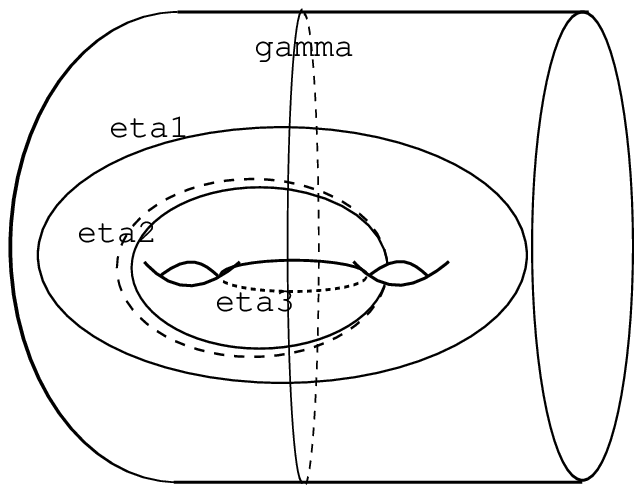} \caption{}
 \end{center}
 \end{figure}

We now explicitly compute $\Psi_2(f_2)$ and see that its characteristic polynomial has degree 2 factors.  Since
$\Psi_2$ is a homomorphism and $\Mod(S)$-equivariant, we find that
$$\begin{array}{rcl} \Psi_2(T_{\gamma_{2,1}} T_{T_P(\gamma_{2,1})}^{-1})& = & \Psi_2(T_{\gamma_{2,1}}) - [T_P]
\circ \Psi_2(T_{\gamma_{2,1}}) \circ
[T_P] ^{-1}\\& = & \Psi_2(T_{\gamma_{2,1}}) - [T_{\eta_1}] [T_{\eta_2}] [T_{\eta_3}]
[T_{\eta_4}]  \Psi_2(T_{\gamma_{2,1}})  [T_{\eta_4}]^{-1} [T_{\eta_3}]^{-1} [T_{\eta_2}]^{-1} [T_{\eta_1}] ^{-1} \\
& = & \Psi_2(T_{\gamma_{2,1}}) - [T_{\eta_1}] [T_{\eta_3}]
 \Psi_2(T_{\gamma_{2,1}}) [T_{\eta_3}]^{-1}  [T_{\eta_1}] ^{-1} \end{array}
$$  Note that since $\eta_2$ is separating, $[T_{\eta_2}]$ is trivial.  For any simple closed curve $\beta$
and $c \in H$, one can show that $$[T_\beta](c) = c + \hat{\iota}([\beta], c) [\beta]$$ where $[\beta]$ is the
homology class of $\beta$.  We see that $[\eta_1] = a_1 + a_2$ and $[\eta_3] = b_2 - b_1$ and so one computes
$$ \Psi_2(T_{\gamma_{2,1}} T_{T_P(\gamma_{2,1})}^{-1}) = 3 *  \left( \begin{array}{rrrr}
    -1 & 0 & 1 & 1 \\
    0 & -1 & 1 & -1 \\
    -1 & -1 & 1 & 0 \\
    -1 & 1 & 0 & 1 \end{array} \right)
    $$
The characteristic polynomial is computed to be $(9 + x^2)^2$ \\

\textbf{Case $g \geq 3$:} First, we find a pair of filling curves using Lemma \ref{lemma:flp}.  Let $\gamma_{g,1}$ be the curve
depicted in Figure 5 and $P$ the pants decomposition depicted in Figure 6.  One sees that $\gamma_{g,1}$ intersects
every curve of $P$ nontrivially.  Thus, by the lemma, $\gamma_{g,1}$ and $\gamma_{g,2} := T_P(\gamma_{g,1})$ fill $S_{g,1}$.

 \begin{figure}[htbp]
 \begin{center}
 \psfrag{gamma}{$\gamma$}
 \psfrag{nu}{$\nu$}
 \psfrag{...}{$\dots$}

 \epsffile{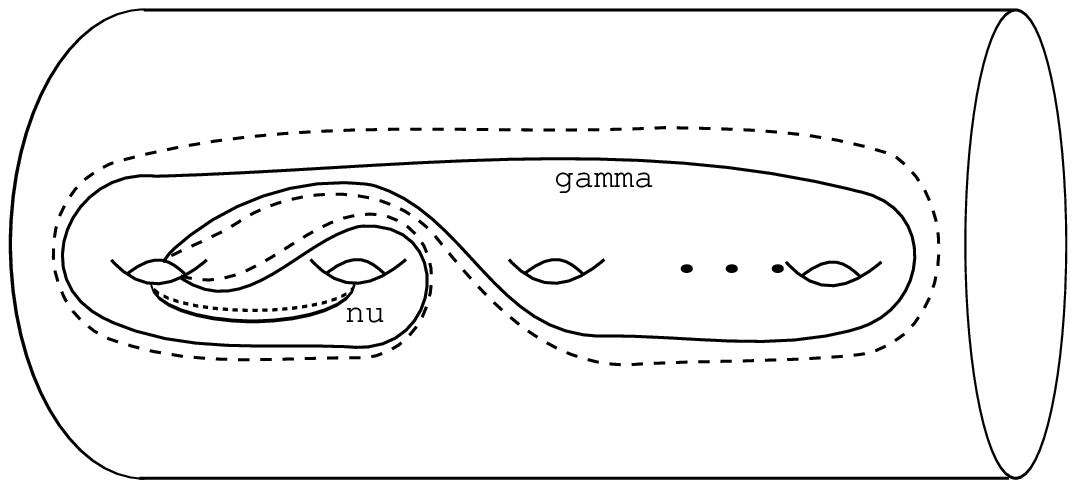} \caption{}
 \end{center}
 \end{figure}

 \begin{figure}[htbp]
 \begin{center}
 \psfrag{...}{$\dots$}

 \includegraphics{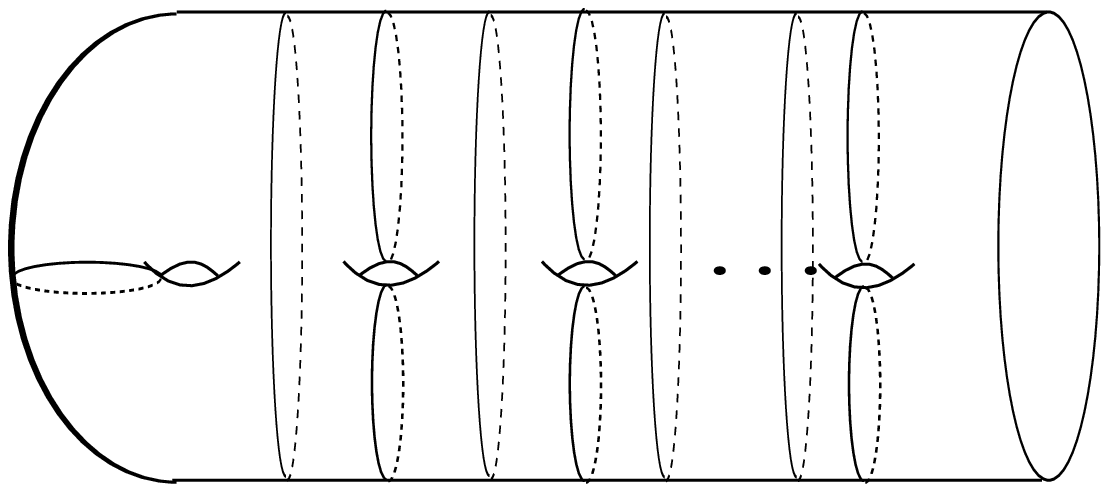} \caption{}
 \end{center}
 \end{figure}

Now, we show that there is some $f_g' \in \mcI_2$ such that $\Psi_2(f_g') = \Psi(f_g)$ and $f_g'$ is reducible.  Let
$$P_{\textrm{nosep}} = \{\eta \in P \;\; | \;\; \eta \; \textrm{is nonseparating}\}$$  Since $[T_\eta] = Id$ for all
$\eta$ that are separating, we see that $$ \begin{array}{rcl} \Psi_2(f_g)& = & \Psi_2(T_{\gamma_{g,1}} T_{T_P(\gamma_{g,1})}^{-1})
= \Psi_2(T_{\gamma_{g,1}}) - [T_P] \Psi_2( T_{\gamma_{g,1}}) [T_P^{-1}]\\ & =& \Psi_2(T_{\gamma_{g,1}}) - [T_{P_{\textrm{nosep}}}] \Psi_2(
T_{\gamma_{g,1}}) [T_{P_{\textrm{nosep}}}^{-1}] = \Psi_2( T_{\gamma_{g,1}} T_{T_{P_{\textrm{nosep}}}(\gamma_{g,1})}^{-1}) \end{array}
$$  We let $f_g' = T_{\gamma_{g,1}} T_{T_{P_{\textrm{nosep}}}(\gamma_{g,1})}^{-1}$.  Notice that the curve $\nu$ in Figure 5
intersects neither
$\gamma_{g,1}$ nor $T_{P_{\textrm{nosep}}}(\gamma_{g,1})$, and so $f_g'(\nu) = \nu$.  Thus, $f'$ is reducible and we are done.
\end{proof}

\subsection{Positive Results for Theorem \ref{theorem:maintheorem}}
In this section, we will exhibit two examples of mapping classes which satisfy the hypothesis of Theorem
\ref{theorem:maintheorem}.  We begin with an example satisfying both Theorem \ref{theorem:maintheorem} and the
Thurston--Penner criteria.

We first make some preliminary remarks.  If $A$ and $B$ are multicurves and $T_A T_B^{-1}$ is pseudo-Anosov, then it
is clear that $A \cup B$ fills $S$.  Thus, if $T_A T_B^{-1}$ satisfies the hypothesis of Theorem
\ref{theorem:maintheorem}, it immediately follows that $T_A T_B^{-1}$ must satisfy the Thurston--Penner criteria.

Now let us describe our example explicitly.  Let $S = S_{5, 1}$.  We let $A = \{\gamma_1, \gamma_2, \gamma_3\}$ and
$B' = \{\gamma_1, \gamma_2\}$ where the $\gamma_i$ are the ``standard'' separating curves given in Figure 3.  Let $h
\in \Mod(S)$ be any mapping class such that its projection to $\Sp(2g, \Z)$ is given by
$$ [h] = \left( \begin{array}{rrrrrrrrrr}
2 & 0 & 1 & 1 & \phantom{0}1 & 1 & \phantom{0}2 & 0 & 1 &  0 \\
1 & 2 & -2& 0 & 0 & -1& 1 & -1& 2 & -2 \\
3 & 3 & 2 & -1& 2 & 0 & 0 & 1 & 2 & -3 \\
1 & -1& 0 & 2 & 1 & 0 & 2 & 0 & 1 & 1 \\
4 & 3 & 2 & -1& 2 & 1 & 1 & 0 & 2 & -2 \\
0 & -1& 2 & 0 & 0 & 1 & 0 & 1 & -1& 1 \\
0 & -1& 0 & 1 & 0 & 0 & 1 & 0 & 0 & 1 \\
6 & 0 & 7 & 2 & 5 & 2 & 3 & 4 & 2 & 0 \\
1 & -1& 2 & 0 & 0 & 1 & 0 & 1 & 0 & 1 \\
1 & 0 & 1 & 0 & 0 & 0 & 0 & 0 & 0 & 1 \end{array} \right)$$ Let $B$ = $h(B')$.  If we let $e_{i, j}$ be the
elementary matrix with a 1 in the $(i, j)$th entry and 0's everywhere else, then using Lemma \ref{lemma:psisepcurve},
we find $$T_A = 15 (e_{1, 1} + e_{2, 2}) + 12 (e_{3, 3} + e_{4, 4}) + 7 (e_{5, 5} + e_{6, 6})$$ and
$$T_{B'} = 8 (e_{1, 1} + e_{2, 2}) + 5 (e_{3, 3} + e_{4, 4})$$
Putting this together, we compute (via Mathematica)
$$ \begin{array}{rcl} \Psi_2 (T_A T_B^{-1})& = & \Psi(T_A) - [h] \Psi(T_{B'}) [h]^{-1} \\
 & = & \left(
 \begin{array}{rrrrrrrrrr}
 42 & 0 & -6 & -33 & -26 & -33 & -25 & 11 & -5 & 26 \\
 0 & 42 & -44 & 14 & -8 & 30 & -116 & 18 & -16 & 24 \\
 14 & 33 & -28 & 0 & -14 & 24 & -89 & 14 & -19 & 38 \\
 44 & -6 & 0 & -28 & -28 & -36 & -22 & 8 & -2 & 20 \\
 30 & 33 & -36 & -24 & -22 & 0 & -89 & 22 & -19 & 46 \\
 8 & - 26 & 28 & -14 & 0 & -22 & 68 & -10 & 8 & -8 \\
 18 & -11 & 8 & -14 & -10 & -22 & 13 & 0 & 3 & 2 \\
 116 & -25 & 22 & -89 & -68 & -89 & 0 & 13 & -10 & 68 \\
 24 & -26 & 20 & -38 & -8 & -46 & 68 & -2 & 8 & 0 \\
 16 & -5 & 2 & -19 & -8 & -19 & 10 & 3 & 0 & 8 \\ \end{array} \right) \end{array} $$
 We compute (via Mathematica) the characteristic polynomial to be $$(x^5 - 21x^4 + 107x^3 +3837x^2 - 13500x + 151200)^2$$
 We find, using Mathematica, that modulo 17 the polynomial $$x^5 - 21x^4 + 107x^3 +3837x^2 - 13500x + 151200$$ is irreducible, and
 hence irreducible over $\Z$.  Thus, by Theorem \ref{theorem:maintheorem}, $T_A T_B^{-1}$ is pseudo-Anosov and we are done.

We now exhibit a mapping class $f \in \mcI_1(S_{4, 1})$ for which there is no obvious way to apply the
Thurston--Penner criteria.  First, let us recall some facts about the Johnson homomorphism on $\mcI_1$.
 There is the following sequence of canonical embeddings and isomorphisms: $$ \Lambda^3 H
\hookrightarrow \Lambda^2 H \otimes H \cong (\Gamma_2 / \Gamma_3) \otimes H \cong \Hom(H, \Gamma_2 / \Gamma_3) $$
Theorem 1 of \cite{J} tells us that
$$\tau(\mcI_1 / \mcI_2) = \textrm{image}(\Lambda^3 H) \subseteq \Hom(H, \Gamma_2 / \Gamma_3)$$
We define a {\em bounding pair} to be a pair of nonisotopic disjoint curves whose union separates the surface.
The {\em bounding pair map} associated to an ordered bounding pair $(\eta, \gamma)$ is the product of Dehn twists
$T_\eta T_\gamma^{-1}$.  Let
 $h = T_{\beta_i} T_{\beta_i'}^{-1}$ be the bounding pair map for $\beta_i$ and $\beta_i'$ as given in Figure 7.

 \begin{figure}[htbp]
 \begin{center}
 \psfrag{beta1}{$\beta_1$}
 \psfrag{beta2}{$\beta_2$}
 \psfrag{betai}{$\beta_i$}
 \psfrag{betag}{$\beta_g$}
 \psfrag{betai'}{$\beta_i'$}
 \psfrag{alpha1}{$\alpha_1$}
 \psfrag{alpha2}{$\alpha_2$}
 \psfrag{alphai}{$\alpha_i$}
 \psfrag{alphag}{$\alpha_g$}

 \includegraphics{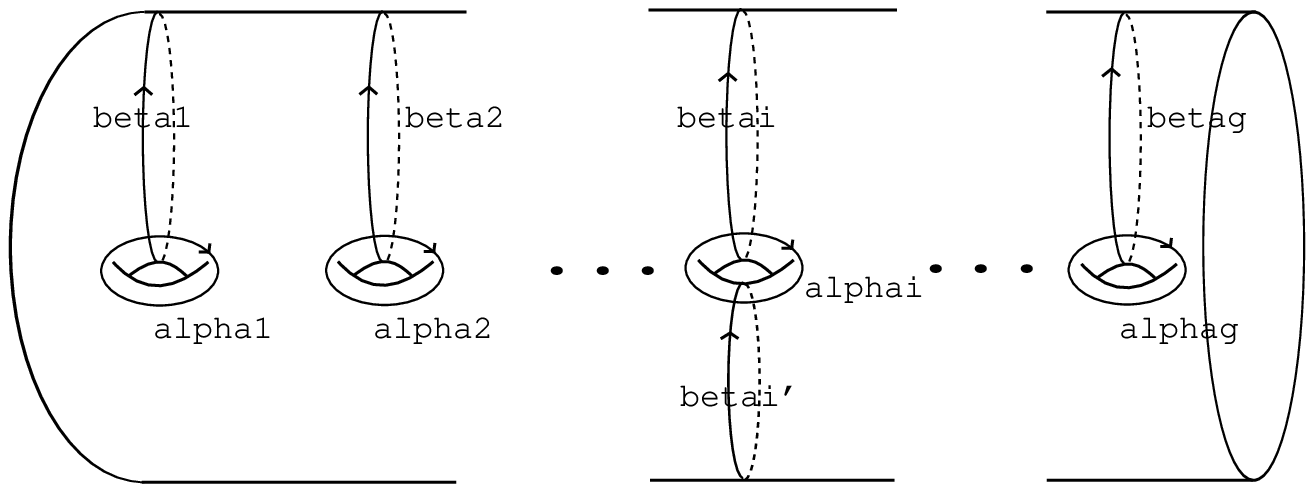} \caption{}
 \end{center}
 \end{figure}

 In Lemma 4B of \cite{J}, Johnson computes that
\begin{equation} \label{eq:Johnson} \tau_h = \left( \sum_{j=1}^{i-1} a_j \wedge b_j \right) \wedge b_i \end{equation}

Now, let us describe the example.  Let
$$\begin{array}{rcl} y & = & (a_4 + b_2 + b_3) \wedge a_1 \wedge b_1 + (a_3 + b_4) \wedge a_2 \wedge b_2 \\  &  & +
(a_1 + a_2 + b_1) \wedge a_3 \wedge b_3 + (a_1 + a_2) \wedge a_4 \wedge b_4
\end{array} \in \Lambda^2 H
$$ From the previous paragraph, we know there exists $f \in \mcI$ such that $\tau_f = y$ which we construct now.
Consider the bounding pairs illustrated in Figures 8.a - 8.h. Let $f$ be the product of bounding pair maps about
these bounding pairs.  Since $\tau$ is a homomorphism to an abelian group, $\tau_f$ is the same regardless of how the
bounding pair maps are composed.  Using (\ref{eq:Johnson}), one computes  that $\tau_f = y$.

 \begin{figure}[htbp]
 \begin{center}
 \psfrag{a}{$a$} \psfrag{b}{$b$} \psfrag{c}{$c$} \psfrag{d}{$d$} \psfrag{e}{$e$} \psfrag{f}{$f$} \psfrag{g}{$g$} \psfrag{h}{$h$}
 \psfrag{1}{1} \psfrag{2}{2}

 \includegraphics{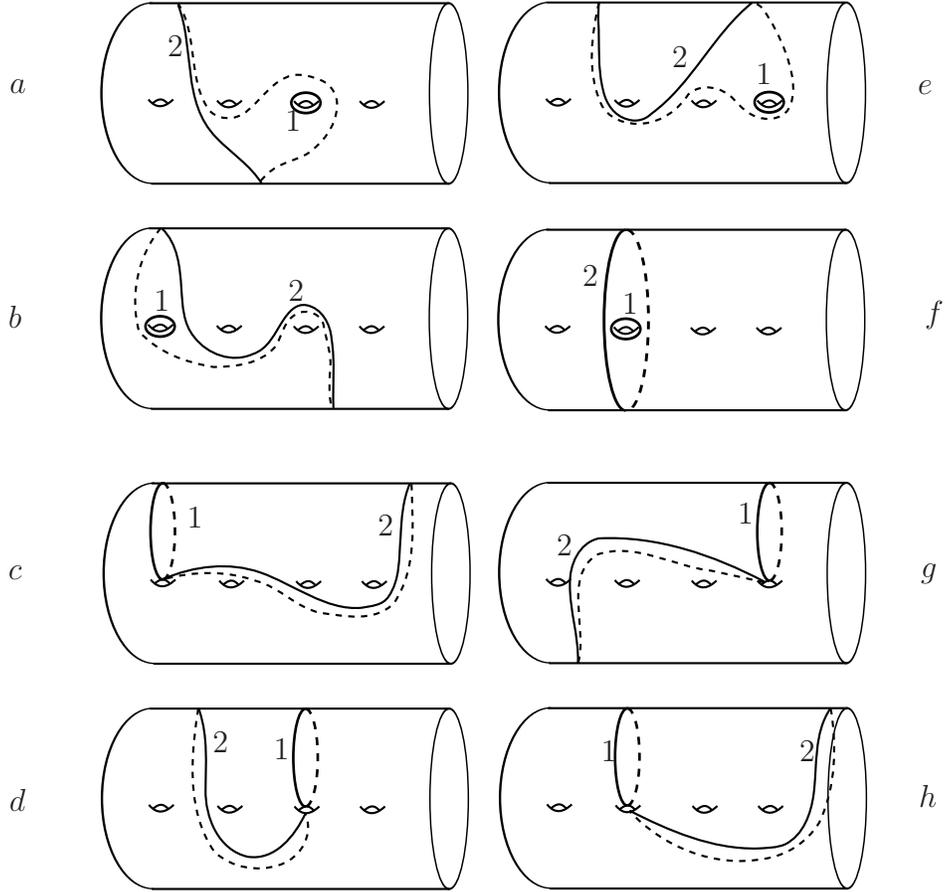} \caption{The product of the bounding pair maps
 indicated in a - h yields $y$ }
 \end{center}
 \end{figure}

 Via computation (with Mathematica),
we find that with respect to the symplectic basis $\{a_1, b_1, \dots, a_4, b_4\}$
$$\Psi_1(f) = \left( \begin{array}{rrrrrrrr} -6 & -2 & 2 & 0 & 2 & 2 & -2 & 0 \\
                                                        4 & 2 & 2 & -2 & -2 & 2 & 2 & -2 \\
                                                        4 & -2 & 2 & 0 & -2 & 2 & 2 & 0 \\
                                                        -2 & 4 & -2 & 0 & 0 & 2 & 0 & 2 \\
                                                        -4 & -4 & 2 & 4 & -2 & 4 & 0 & -2 \\
                                                        -4 & -4 & 0 & 6 & 2 & 2 & -2 & 0 \\
                                                        -2 & 4 & -2 & 2 & 2 & 2 & 2 & 2 \\
                                                        4 & -2 & -2 & -4 & -4 & 2 & 4 & 0 \end{array}
                                                        \right)  $$
The characteristic polynomial of $\Psi_1(f)/2$ is $$\chi(\Psi_1(f)/2) = x^8 - 8x^6 + 26x^5 - 18x^4 -76x^3 +241x^2
-558x +553$$ This polynomial is found to be irreducible mod 11 via Mathematica and is hence irreducible.  By Theorem
\ref{theorem:maintheorem}, $f$ is pseudo-Anosov. Note that curves $c_2$, $d_2$, and $g_2$ in Figure 8 all pairwise
intersect, and so the criteria of Thurston and Penner do not seem to apply directly to $f$.

\end{document}